\documentclass[12pt, twoside, leqno]{article}

\usepackage{amsmath,amsthm}
\usepackage{amssymb}

\usepackage{enumerate}

\usepackage{graphicx}

\usepackage{pgf,tikz}
\usetikzlibrary{arrows}

\newtheorem{thm}{Theorem}[section]

\newtheorem{prop}[thm]{Proposition}

\theoremstyle{definition}
\newtheorem{defin}[thm]{Definition}
\newtheorem{rem}[thm]{Remark}
\newtheorem{exa}[thm]{Example}

\numberwithin{equation}{section}

\usepackage{pgf,tikz}
\usetikzlibrary{arrows}

\begin{document}

\baselineskip=12.1pt


\title{On the   distance eigenvalues of design graphs}

\author{S. Morteza Mirafzal\\
Department of Mathematics \\
  Lorestan University, Khorramabad, Iran\\
E-mail: mirafzal.m@lu.ac.ir\\
E-mail: smortezamirafzal@yahoo.com}

\date{}

\maketitle

\renewcommand{\thefootnote}{}


\footnote{\emph{Date}:  }

\renewcommand{\thefootnote}{\arabic{footnote}}
\setcounter{footnote}{0}
\date{}

\begin{abstract} 
 A design graph is a regular bipartite graph in which any two distinct
vertices of the same part  have the same number of common neighbors. This class of graphs have a close relationship to strongly regular graphs. In this paper, we    study   the  distance eigenvalues of the design graphs. Also, we will explicitly determine  the  distance eigenvalues of a  class of design graphs, and determine the values for which the class is   distance integral, that is,  its distance eigenvalues are integers.

\end{abstract}

\maketitle

\section{ Introduction and Preliminaries}
  In this paper, a graph $G=(V,E)$ is
considered as an undirected simple graph where $V=V(G)$ is the vertex-set
and $E=E(G)$ is the edge-set. For all the terminology and notation
not defined here, we follow [4,6,7,8,20].\

Let $G=(V,E)$ be a graph and  $A=A(G)$ be an adjacency matrix of $G$. The $characteristic$  $polynomial$  of $G$ is defined as $P(G; x)=P(x) = det(xI-A)$. A zero of $P(x)$ is called an $eigenvalue$  of the graph $G$. A graph
is called $integral$  if all its eigenvalues are integers. The study of integral graphs  was initiated by Harary
and Schwenk in 1974 (see [9]). A survey of papers up to 2002 has been  appeared in [3].  \newline 
Let $n$ be the number of vertices of the graph $G$. The distance matrix $D=D(G)$ is an $n \times n $ matrix indexed by $V$  such that $D_{uv} = d_{G}(u, v)=d(u,v)$  where $d_{G}(u, v)$ is the distance between the vertices $u$ and $v$ in the graph $G$.
The characteristic polynomial $P(D; x) = det(xI - D)=D_G(x)$ is the $distance$ $characteristic$  $polynomial$  of $G$. Since $D$ is a real
symmetric matrix, the distance characteristic polynomial $D_G(x)$ has real zeros. Every zero of the polynomial $D_G(x)$ is called a $distance$ $eigenvalue$  of the graph $G$. 
Two surveys on the distance
spectra of graphs are [1,10].  An interested reader may see some of the recent works on  distance spectra of graphs  in [2,5,11,18,20]. \

A graph $G$ is called $distance$ $integral$ (briefly, $D$-integral) if all the distance eigenvalues  of $G$ are integers. The $D$-integral graphs are studied only in a few number of papers (see [17,18,20]). \

Let $G=(V,E)$ be a graph with the vertex-set $V=\{ v_1,\dots,v_n \}$ and distance matrix  $D(G)=D=(d_{ij})_{n \times n}$, where $d_{ij}=d(v_i,v_j)$. Let $H \leq Aut(G)$  and 
 $ \pi = \{ w_1^H=C_1,\dots,w_m^H=C_m \} $ be the orbit partition of $H$, where $\{w_1,\dots,w_m  \} \subset V$. Let $Q=Q_{\pi}=(q_{ij})_{m \times m}$ be the matrix which its rows and columns are indexed by $\pi$ such that,
   $$q_{ij}=\sum _{w\in C_j} d(v,w), \ \ \ \ \ (*) $$
     where $v$ is a fixed element in the cell $C_i$. It is easy to check that this sum is independent of $v$, that is, if $u \in C_i$, then $q_{ij}=\sum _{w\in C_j} d(v,w)=\sum _{w\in C_j} d(u,w)$. Hence, the matrix $Q$ is well defined. We call the matrix $Q$ the $quotient \ matrix $ of $ D $  over      $\pi  $.  We claim that every eigenvalue of $Q$ is an eigenvalue of the distance matrix $D$. In fact we have the following facts. 
 
  \begin{thm}  $[17]$  Let $G=(V,E)$ be a graph with the distance matrix $D$. Let  $\pi$ be an orbit partition of $V$    and $Q$ be a quotient matrix of $D$ over $\pi$.  Then  every eigenvalue of $Q$ is an eigenvalue of the distance matrix $D$.
  
\end{thm}
   
\begin{prop}$[17]$ Let $G=(V,E)$ be a graph with the distance matrix $D$. Let  $\pi$ be an orbit partition of $V$   and $Q$ be a quotient matrix of $D$ over $\pi$.   Let $\theta$ be an eigenvalue of the distance matrix $D$  with the 
non zero eigenvector $f$ such that $f$ is constant on every cell of $\pi$.    Then  $\theta$ is an eigenvalue of the   matrix $Q$. 

\end{prop}

We now can easyly verify  the following statement by an easy discussion [17].
  
\begin{thm}Let $G=(V,E)$ be a   graph and $D$ be a distance matrix for $G$. Let $f \neq 0$ be an eigenvector with the eigenvalue $\lambda$ for $D$. Let $H$ be a subgroup of $Aut(G)$ and $\pi$ be its orbit partition on $V$ and $Q$ is a quotient matrix of $D$ over $\pi$. If $\lambda$ is not an eigenvalue of $Q$, then the sum of 
the values of $f$ on each cell of $\pi$ is zero.

\end{thm}
   The following fact is an important tool for finding some interesting results in the present wok.  
\begin{thm} 
  Let $G=(V,E)$ be a vertex-transitive graph with the distance matrix $D$.
  Let $H$ be a subgroup of $Aut(G)$ with the orbit partition $\pi$ on $V$ such that $\pi$ has a singleton cell $\{x\}$. Let $Q=Q_{\pi}$ be a quotient   matrix of $D$ over $\pi$. Then the set of distinct eigenvalues of $D$  is  equal to the set of distinct eigenvalues of $Q$. 
\end{thm}

   A partition $V=V_1 \cup V_2 \cup \dots \cup V_m$ is called a $distance$ $equitable$ $partition$ for the vertex-set of the graph $G=(V,E)$    if for every $v \in V_i$ the sum $\sum _{w \in V_j} d(v,w)$
is constant, that is, $\sum _{w \in V_j} d(v,w)$ is independent of the choice of the vertex $v$. Hence, as we saw, every orbit partition for the vertex-set $V$ is a distance equitable partition. Now if we define by the partition $\Pi=\{ V_1, \dots, V_m \}$ the matrix $P=P{_\Pi}=(p_{ij})_{m \times m}$ by the rule, $$p_{ij}=\sum _{w\in V_j} d(v,w), $$
where $v$ is a fixed element in the cell $C_i$, then  the matrix $P$ is well defined. 
We call $P$ a quotient distance matrix over the partition $\Pi$, or a quotient matrix for the distance matrix $D$ over $\Pi$. 
By a similar argument which have been appeared in the proof of Theorem 2.1,   in [17],  we can check that the following fact holds.

 \begin{thm}    Let $G=(V,E)$ be a graph with the distance matrix $D$. Let  $\Pi$ be a distance equitable partition of $V$  and $P$ be a quotient matrix of $D$ over $\Pi$.  Then   every eigenvalue of $P$ is an eigenvalue of the distance matrix $D$.
  
\end{thm}
We have a result similar to what is appeared in Theorem 1.1, for the case distance  equitable partitions.

\begin{prop}  Let $G=(V,E)$ be a graph with the distance matrix $D$. Let  $\Pi$ be a distance equitable partition of $V$   and $P$ be a quotient matrix of $D$ over $\Pi$.   Let $\theta$ be an eigenvalue of the distance matrix $D$  with the 
non zero eigenvector $f$ such that $f$ is constant on every cell of $\Pi$.    Then  $\theta$ is an eigenvalue of the   matrix $P$. 

\end{prop} 

\begin{proof}
The proof is similar to the discussion which has been appeared in [17, Prop 2.2].
\end{proof}

We now introduce a class of graphs,  one which some classes of graphs with interesting algebraic properties  are subclasses of it.

\begin{defin}$[19]$
A design graph with parameters $(m,d,c)$, $c \neq 0$,  is a $d$-regular bipartite graph of order $2m$ in which any two distinct
vertices of the same part have  $c$  common neighbor(s). The
complete bipartite graph $K_{n,n}$ fits the definition, but we exclude it by convention.
\end{defin}
It is easy to check (by double counting method) that  if $G=(V,E)$ is a design graph with parameters $(m,d,c)$, then  the following equality holds, 
$$c(m-1)=d(d-1).$$
It follows from Definition 1.7, that  every design graph is a connected graph with diameter  3 and girth 4 or 6 according to whether  $c>1$ or $c=1$.  Conversely, a regular bipartite graph of diameter 3 and
girth 6 is a design graph with parameter $c = 1$ [19].\

 \begin{exa}Let $n\geq 3$   be  an integer. Let  $V$ be the set of all $1$-subsets
and $(n-1)$-subsets of $[n]=\{1,\dots,n \}$.
 The $bipartite\  Kneser\  graph$ $H(n, 1)$ has
$V$ as its vertex-set, and two vertices $v,  w$ are adjacent if and only if $v \subset w$ or $w\subset v$. It is easy to see that $H(n, 1)$ is a bipartite graph of diameter 3. Some of properties of the graph $H(n,1)$ and a generalization of it have been appeared in [13,14,15,16].
We can easyly show that $H(n,1)$ is a design graph with parameters $(n,n-1,n-2)$.

 \end{exa}

 A  $strongly \ regular$ graph  with parameters $(n,k,a,c)$  is a  $k$-regular graph of order $n$ in which every pair of adjacent vertices has exactly $a$ common neighbors and every pair of nonadjacent
vertices has exactly $c$ common neighbors. For example,  the cycle $C_5$ is a strongly regular graph with parameters $(5,2,0,1)$ and the Petersen graph is a strongly regular graph with parameters $(10,3,0,1)$.\

Let $G=(V,E)$ be a graph. The $bipartite\  double\  cover$ of $G$ which we denote it  by $B(G)$  is a graph with the vertex-set $V \times \{ 0,1 \}$,     in which vertices $(v,a)$ and $(w,b)$ are adjacent if and only if 
 $a \neq b$ and $\{v,w\} \in E$.  For example, if $n >2$ is an odd integer, then  the bipartite double cover of the cycle $C_n$ is (isomorphic with) the cycle $C_{2n}$. Also, the bipartite double cover of the complete graph $K_n$ is (isomorphic with) the bipartite Kneser graph $H(n,1)$.
\begin{exa}
Let $G=(V,E)$ be a strongly regular graph with parameters $(n,d,a,c)$ with $a\neq 0$ and $a=c$. It is not hard to check that the bipartite double cover of $G$, that is, the graph $B(G)$ is a design graph with parameters $(n,d,c)$.

\end{exa}
In this paper, we will study   the distance eigenvalues of design graphs. Also we will explicitly determine the distance eigenvalues of an important class of design graphs and determine the values for which the class is integral. The main tool which we use in our work is the distance equitable partition and the  orbit partition method in
algebraic graph theory,  which we have already employed it in determining the adjacency eigenvalues of a particular family of  graphs [12]  and later in determining the sets of distance eigenvalues of some other families of graphs [17,18].    We will show how we can find, by using this method, the set
of all distinct distance eigenvalues of  some classes of design graphs.

\section{ Main Results}
Let $G=(V=V_1 \cup V_2,E)$, $V_1 \cap V_2=\emptyset$ be a  design graph with parameters $(m,d,c)$,   the adjacency matrix $A=(a_{ij})$ and distance matrix $D=(d_{ij})$.  Consider the partition $\pi=\{ V_1,V_2\}$. If $v \in V_1$, then the number of its neighbors in $V_2$ is $d$.  If $w \in V_2$ is not adjacent to $v$ then it is at distance 3 from $v$. Thus  we have,  \\\\
$\sum_{w \in V_2}d(v,w)=d+3(m-d)=3m-2d.$\\\\
On the other hand, if $w \in V_1$, since $v$ and $w$ have at list a common neighbor, then $d(v,w)=2$. Hence we have, \\\\
 $\sum_{w \in V_1}d(v,w)=2(m-1).$\\\\
We now deduce that $\pi=\{ V_1,V_2\}$ is a distance  equitable  partition for the vertex-set $V.$  Hence, the following matrix $Q$ is a quotient matrix of the distance matrix $D$ over the partition $\pi$,  
$$Q=\begin{pmatrix}
2(m-1)&3m-2d\\
3m-2d&2(m-1)
\end{pmatrix}.$$

It is easy to see that each of the functions (vectors) $f_1=(1,1)^t$ and  $f_2=(1,-1)^t$ are the eigenvectors for the matrix $Q$ with the eigenvalues,  \\\\ $\gamma_1=2(m-1)+3m-2d=5m-2d-2$ and \\\\ $\gamma_2=2(m-1)-(3m-2d)=-m+2d-2$, \\\\ respectively. Hence,    from Theorem 1.5,   we have the following result.

\begin{prop} Let $G=(V_1 \cup V_2,E)$, $V_1 \cap V_2=\emptyset$ be a  design graph with parameters $(m,d,c)$. Then,  the integers  ${\gamma}_1=5m-2d-2$ and ${\gamma}_2=-m+2d-2$ are distance eigenvalues of $G$.  

\end{prop}

\begin{rem}
Concerning Proposition 2.1, we can easyly check  that the eigenvector corresponding to the eigenvalue $\gamma_1$ for the distance matrix $D$ of the graph $G$  is the function $e_1$, defined by the rule $e_1(v)=1$, for every $v \in V_1 \cup V_2=V(G)$. Also, the eigenvector corresponding  to the eigenvalue $\gamma_2$  for the distance matrix $D$ of the graph $G$  is the function $e_2$, defined by the rule $e_2(v)=1$, for every $v \in V_1$, and $e_2(v)=-1$, for every $v \in V_2.$ 

\end{rem}

We now proceed to construct another distance  equitable   partition for the vertex-set of the design graph $G=(V,E), \  V=V_1 \cup V_2, \ V_1 \cap V_2=\emptyset,$ with parameters $(m,d,c).$ 
Let $v_1$ be an arbitrary vertex in the part $V_1.$ We let, \\\\ $O_1=\{ v_1 \}$, \\\\ $O_2=V_1-\{v_1\}$,  \\\\$O_3=N(v_1)$, the set of neighbors of $v$ in $V_2, $   and \\\\ $O_4=V_2-N(v_1)=V_2-O_3.$ \\\\ We claim that  $\pi_2=\{ O_1,O_2,O_3,O_4 \}, $   is a distance equitable partition for the vertex-set $V$ of $G$. Let $v_i \in O_i, $ $i \in \{2,3,4 \}$ be arbitrary vertices. Then we have the following equalities.\\\\
$p_{11}=\sum_{w \in O_1}d(v_1,w)=0.$\\\\
$p_{12}=\sum_{w \in O_2}d(v_1,w)=2(m-1).$\\\\
$p_{13}=\sum_{w \in O_3}d(v_1,w)=d.$\\\\
$p_{14}=\sum_{w \in O_4}d(v_1,w)=3(m-d).$\\\\
$p_{21}=\sum_{w \in O_1}d(v_2,w)=2.$\\\\
$p_{22}=\sum_{w \in O_2}d(v_2,w)=2(m-2).$\\\\
$p_{23}=\sum_{w \in O_3}d(v_2,w)=c+3(d-c)=3d-2c.$\\
 Because vertices $v_1$ and $v_2$ have $c$ common neighbors which are in $N(v_1)=O_3.$ Moreover, if the vertex $w$ in $O_3$ is not adjacent to $v_2, $ then $d(v_2,w)=3.$\\\\
$p_{24}=\sum_{w \in O_4}d(v_2,w)=(d-c)+3((m-d)-(d-c))=3(m-d)-2d+2c=3m-5d+2c.$\\
 Because $c$ neighbors of  $v_2$ are in $O_3$, and hence $d-c$ neighbors of it are in $O_4.$ \\\\
$p_{31}=\sum_{w \in O_1}d(v_3,w)=1.$\\\\
$p_{32}=\sum_{w \in O_2}d(v_3,w)=(d-1)+3((m-1)-(d-1))=3m-2d-1.$\\
Because, $(d-1)$ neighbors of $v_3$ are in $O_2.$\\\\
$p_{33}=\sum_{w \in O_3}d(v_3,w)=2(d-1).$\\\\
$p_{34}=\sum_{w \in O_4}d(v_3,w)=2(m-d).$\\\\
$p_{41}=\sum_{w \in O_1}d(v_4,w)=3.$\\\\
$p_{42}=\sum_{w \in O_2}d(v_4,w)=d+3(m-1-d)=3m-2d-3.$\\\\
$p_{43}=\sum_{w \in O_3}d(v_4,w)=2d.$\\\\
$p_{44}=\sum_{w \in O_4}d(v_4,w)=2(m-d-1).$\\\\
Now from our argument it follows that $\pi_2$ is a distance equitable partition for the vertex-set $V.$
Hence,  we have the following result.

\begin{thm} Let $G=(V,E), \ V=V_1 \cup V_2, \ V_1 \cap V_2 =\emptyset$, be a design graph with parameters $(m,d,c)$. Let $v$ be a vertex in $V_1$. Then the partition
 $$\pi_2=\{O_1=\{v\}, O_2=V_1-O_1, O_3=N(v), O_4=V_2-O_3\}, \ \ \ \ (*)$$
  is a distance equitable partition for the vertex-set $V$ and the following matrix $P=(p_{ij})_{4 \times 4}$ is a distance quotient matrix over $\pi_2$, 
  
$$P=P(m,d,c)=\begin{pmatrix}
0&2m-2&d&3m-3d\\
2&2m-4&3d-2c&3m-5d+2c\\
1&3m-2d-1&2d-2&2m-2d\\
3&3m-2d-3&2d&2m-2d-2

\end{pmatrix}. \ \ \ \ (**)$$

\end{thm}

\begin{rem} From Proposition 2.1, remark 2.2, and proposition 1.6, we can deduce that  ${\gamma}_1=5m-2d-2$ and ${\gamma}_2=-m+2d-2$ are distance eigenvalues of the matrix $P$ which is defined in Theorem 2.3. Now one can find the other eigenvalues of $P$ by some handy calculations.

\end{rem}

We now conclude  from Theorem 1.4,   Theorem 1.5, and Theorem 2.3,   the following important result.
\begin{thm}Let $G=(V,E), \ V=V_1 \cup V_2, \ V_1 \cap V_2 =\emptyset$, be a design graph with parameters $(m,d,c)$. Then, 
every eigenvalues of the matrix $P$ defined in Theorem $2.3$, is a distance eigenvalue for the graph $G$. Moreover, if the design graph $G$ is vertex-transitive and  the partition $\pi_2$ defined in Theorem $2.3$, is an orbit partition,   then the set of   eigenvalues of $P$ is   equal  to the set of   distance eigenvalues of the graph $G$. 

\end{thm}

In the next section, we will see how Theorem 2.5, can help us to determine the distance eigenvalues of a design graph. 

\section{An Application}

Let $q$ be a power of a prime $p$ and $\mathbb{F}_q$ be a finite field of order $q$. Let $V(q,n)$ be a vector space of dimension $n$ over $\mathbb{F}_q$. We define  the graph $S(q,n,1)$ as  a  graph with the vertex-set $V=V_1 \cup V_{n-1}$, where $V_1$ and $V_{n-1}$ are the family of subspaces in $V(q,n)$ of dimension $1$ and $n-1, $ respectively,     in which two vertices $v$ and $w$ are adjacent whenever $v$ is a subspace of $w$  or $w$ is a subspace of $v$. 
    It is clear that $S(q,n,1)$  is a bipartite graph with partition  $V=V_1\cup V_{n-1}$. 

We know that the number of $k$-subspaces of a vector space $U$ of dimension $n$ over the  filed $\mathbb{F}_q$ is the Gaussion number   
$${n\brack k}_q= \dfrac{(q^{n}-1)(q^n-q)\dots (q^{n}-q^{k-1})}{(q^{k}-1)(q^k-q)\dots (q^k-q^{k-1})} =\dfrac{(q^{n}-1)\dots (q^{n-k+1}-1)}{(q^{k}-1)\dots (q-1)}. $$
Thus, $|V_1|={n \brack 1}$ and $|V_{n-1}|={n\brack n-1}={n \brack 1} $.\\\\
Hence,  the graph $S(q,n,1)$ is of order $2{n \brack 1}$. It can be checked that $S(q,n,1)$ is a connected regular graph of diameter 3. It is not hard to show that $S(q,n,1)$ is a vertex-transitive graph [16,19]. Also, it can be check that the girth of the graph $S(q,n,1)$ is 6 when $n=3$ and is 4 when $n \geq 4$.  Some properties of $S(q,3,1)$ have been investigated in [8, Chaper 5].  It is easy to show that the graph $S(q,n,1)$ is a design graph with parameters $(m,d,c)$ [19],  where
 $$m=\frac{q^n-1}{q-1}, \  d=\frac{q^{n-1}-1}{q-1}, \  c=\frac{q^{n-2}-1}{q-1}.\ \ \ \ (1)$$ 

Let $K$ be a field,   $V(K,n)$  a vector space of dimension $n$ over and $K$ and $GL(n,K)$ be the group of non-singular linear mappings on the space  $V(K,n)$. 
  This group  contains a normal subgroup isomorphic to $K^{*}$,
  namely,  the group $Z= \{ kI_{V(K,n)} | k \in K \}$, where $I_{V(K,n)}$ is the identity mapping on $V(K,n)$. We denote the quotient group  $\frac{GL(n,K)}{Z}$  by 
$PGL(n,K)$.\

Note that if $(a+Z) \in PGL(n,K)$ and $x$ is an $m$-subspace of $V(K,n)$, then $ (a+Z)(x)=\{ a(u) | u \in x \}$ is an $m$-subspace of  $V(K,n)$. In the sequel, we also denote $(a+Z) \in P GL(n,K)$ by $a$. Let $V_m$ denote the set of $m$-subspaces of the vector spaces $V(n,k).$ 
 Now, if $a \in PGL(n,K)$,  it is easy to see that  the mapping 
  $f_a : V_m \longrightarrow V_m$, defined by   the  rule $f_a(v)=a(v)$,  
  is a well defined function.  Therefore  if we let 
   $$A(n)=\{ f_a | a \in PGL(n,K)  \},   \ \ \ \ \ \ (2) $$
    then $A(n)$ is a group isomorphic to the group $PGL(n,K)$ (as abstract groups),  which acts transitively on the set $V_m$. \
    
Now, it is easy to check that $A(n)$ is a subgroup of automorphism group of the graph $S(q,n,1).$   
 The graph $S(q,n,1)$ has some other automorphisms.  In fact,  the mapping $t$ defined on its vertex-set by the rule
    $t(v)=v^{\perp}$, where   $v^{\perp}$ is the orthogonal complement of $v$,   is an automorphism of the graph $S(q,n,1)$. 
    Hence $M=\langle A(n),t \rangle   \leq Aut(G)$. 
     Note that the order of $t$  is 2 and hence $\langle t \rangle   \cong \mathbb{Z}_2$.
      Let $w\in V_1$ be a vertex in the graph $S(q,n,1).$ Let $H$ be a subgroup of $A(n)$ which fixes $w.$   Now, it is not hard to check that the following partition $\pi$  is an orbit partition generated by $H$ on the vertex-set of $S(q,n,1),$
$$\pi=\{O_1=\{w\}, O_2=V_1-O_1, O_3=N(w), O_4=V_2-O_3\}. \ \ \ \ $$     
Since the graph $S(q,n,1)$ is a vertex-transitive graph, it follows from Theorem 2.4, that  the matrix $P=P(m,d,c)$ defined in the theorem contains all the distance eigenvalues of the graph $S(q,n,1), $ where $m,d,c$ are parameters which are defined in (1). 
We can see  from  Remark 2.4, (or Theorem 1.3, and Remark 2.2, because $\pi$ is an orbit partition for the vertex-set of the graph $ (S(n,q,1))) $ that, 
$$\lambda_1=5\frac{q^n-1}{q-1}-2\frac{q^{n-1}-1}{(q-1)}-2=\frac{5q^n-2q^{n-1}-2q-1}{q-1}, $$
and
$$\lambda_2=-\frac{q^n-1}{q-1}+2\frac{q^{n-1}-1}{q-1}-2=\frac{-q^n+2q^{n-1}-2q+1}{q-1}, $$
are distance eigenvalues of the graph $S(q,n,1).$ Having these eigenvalues in the hand, we can find two other eigenvalues of the matrix $P=P(m,d,c)$ by some  handy calculations. We can also find 
the eigenvalues of the matrix $P$ by a suitable software program. Using Wolfram Mathematica [21]  we have, \\\\
$P := \{\{0, 2m-2, d,3 m - 3 d \},\\
 \{2, 2 m - 4, 
   3d-2c,3 m - 5 d + 2 c\}, \\
    \{1,3m- 2 d - 1, 
   2 d - 2 , 2m-2d\}, \\
   \{3,3m-2d-3,2d,2m-2d-2\} \}$\\\\
$m := (q^n - 1)/(q - 1)\\\\
d := (q^{n - 1} - 1)/(q - 1)\\\\
c := (q^{n - 2}- 1)/(q - 1)$\\\\
Eigenvalues[P]=\{
$r_1=- \frac{q + 2 q^2 + 2 q^n - 5 q^{1 + n}}{(-1 + q) q},$
$r_2=- \frac{-q + 2 q^2 - 2 q^n + q^{1 + n}}{(-1 + q) q},$ \\\\\\
$r_3=\frac{2( q^2 - q^3 - \sqrt{q^{2 + n} - 2 q^{3 + n} + q^{4 + n}})}{(-1 + 
   q) q^2},$
   $r_4=\frac{2( q^2 - q^3 + \sqrt{q^{2 + n} - 2 q^{3 + n} + q^{4 + n}})}{(-1 + 
   q) q^2}  \  \}$. \\\\\\
   Note that $r_1=\lambda_1$ and $r_2=\lambda_2$. 
   Hence, the set $R=\{r_1,r_2,r_3,r_4 \}$ is the set of distance eigenvalues of the graph $S(n,q,1).$ Note that $r_1$  and $r_2$ are integers, since they are roots of a monic polynomial with integer coefficients, namely, they are algebraic integers, and every rational algebraic integer is an integer. Also, we have $q^{2 + n} - 2 q^{3 + n} + q^{4 + n}$=$q^{n+2}(1-2q+q^2)=q^{n+2}{(q-1)}^2.$ Hence,\\\\
    $r_3=\frac{2(q^2(1-q)-(q-1)\sqrt{q^{n+2}})}{q^2(-1+q)}$=$\frac{-2q^2-\sqrt{q^{n+2}}}{q^2}$,  and \\\\
    $r_4=\frac{-2q^2+\sqrt{q^{n+2}}}{q^2}.$ \\\\
     Now it follows that if $n$ is an even integer, then 
the graph $S(n,q,1)$ is a distance integral graph. We now have the following result.

\begin{thm} Let $n \geq 3$ be an integer. Then the set of distance eigenvalues of the graph $S(n,q,1)$ is the set,\\
$$\{\frac{5q^n-2q^{n-1}-2q-1}{q-1}, \frac{-q^n+2q^{n-1}-2q+1}{q-1},$$ \\
$$ \frac{-2q^2-\sqrt{q^{n+2}}}{q^2},\frac{-2q^2+\sqrt{q^{n+2}}}{q^2} \}.$$
Hence, if $n$ is an even integer, then the graph $S(n,q,1)$ is a distance integral graph.
\end{thm}

\section{Conclusion}
In this paper, we   studied  the distance eigenvalues of design graphs. In particular, by finding the set of distinct distance eigenvalues of a class of design graphs, we showed that the class is a distance integral graph (Theorem 3.1).  The main tools which we employed were the distance equitable  and  orbit partition methods in algebraic graph theory.

\end{document}